\documentclass[twocolumn,11pt]{article}
\usepackage{times}
\usepackage{graphicx}
\def\NN{{\rm I\hskip-2pt N}}
\def\RR{\vbox {\hbox to 8.9pt {I\hskip-2.1pt R\hfil}}}
\def\CC{{\rm C\hskip-4.8pt \vrule height 6pt width 12000sp\hskip 5pt}}

\def\e{\hbox {e}}
\def\exp{{\rm exp}}
\def\ds{\displaystyle}
\def\rec#1{\frac{1}{#1}}

\def\pni{\par \noindent}
\def\vsh{\vskip 0.25truecm\noindent}

\def\vsp{\vsh\pni}

\begin{document}
\global\def\refname{{\normalsize \it References:}}
\baselineskip 12.5pt
%
%
%
\title{\LARGE \bf The Pioneers of the Mittag-Leffler  Functions \\ in 
Dielectrical and Mechanical Relaxation Processes}

\date{Paper published with footnotes in WSEAS Transactions on Mathematics \\
Vol 19 (2020), pp. 289--300. DOI: 10.37394/23206.2020.19.29}

\author{ 
\begin{minipage}[t]{3.3in} \normalsize \baselineskip 12.5pt
\centerline{~~~~~FRANCESCO MAINARDI}
\centerline{~~~~~~~University Bologna and INFN}
\centerline{~~~~Department of Physics and  Astronomy}
\centerline{~~~~francesco.mainardi@bo.infn.it}
\end{minipage} 
\begin{minipage}[t]{3.3in} \normalsize \baselineskip 12.5pt
	\centerline{ARMANDO CONSIGLIO~~}
	\centerline{Universit\"{a}t W\"{u}rzburg}
	\centerline{Institut f\"{u}r Theoretische Physik
	 und Astrophysik}
	\centerline{~~armando.consiglio@physik.uni-wuerzburg.de}
\end{minipage}
\kern 0in
%
\\ \\ \\   
\begin{minipage}[b]{6.9in} \normalsize
\baselineskip 12.5pt {\it Abstract:}
%
We start with a short survey of the basic properties of the Mittag-Leffler  functions. 
 Then we  focus on the key role of these  functions to explain the after-effects and relaxation phenomena occurring in dielectrics and  in viscoelastic bodies.
 For this purpose  we recall  the main aspects that were  formerly discussed by two pioneers in the years 1930's-1940's whom 
 we have identified with  Harold T. Davis\footnote{{Harold Thayer Davis (1892-1974) was a mathematician
 	interested in Differential and Integral Equations, Special
 	Functions, Statistics, Econometrics.
 	He was Professor of Mathematics at Indiana University,
 	Bloomington and later at Northwestern University, Illinois.}} ~and Bernhard Gross\footnote{{Bernhard Gross (1905-2002) was an interdisciplinary physicist
 	usually known for Linear Systems in Dielectrics and
 	Viscoelasticity. As a young physicist from Stuttgart, in 1933 he
 	established at the Technological Institute in Rio de Janeiro.
 	Later, in 1960's, for seven years he was Director of the Division
 	of Scientific and Technological Information at the International
 	Atomic Energy Agency (IAEA) in Vienna. Then, after returned to
 	Brazil from Vienna, he started working at Sao Carlos, University
 	of Sao Paulo.}}.
\\ [4mm] {\it Key--Words: Mittag-Leffler functions,
Laplace transforms,  Dielectrics, Electrical Circuits, Viscoelasticity, After-effect, Relaxation.}
\end{minipage}
\vspace{-10pt}
}

\maketitle

\thispagestyle{empty} \pagestyle{empty}
%
%
\section{Introduction}
The functions of the Mittag-Leffler type are known to  play a very important role in the applications of Fractional Calculus, so the original Mittag-Leffler function is recalled as {\it the Queen function of the Fractional Calculus} as it was suggested by Mainardi \& Gorenflo in 2007
\cite{Mainardi-Gorenflo FCAA07}.
\\
The Fractional Calculus is the mathematical theory of the generalized calculus
including operators interpreted as integrals and derivatives of non integer order.
It has a long history starting from Leibniz and involving eminent mathematicians like
Abel, Euler, Fourier, Weyl, Liouville, Riemann just to cite the most prominent ones  of the 17-th, 18-th, 19-th  centuries.  The interested reader can consult any book on Fractional Calculus published after 1974, including those
  (in order of publication date up to 2010) by:
Oldhan \& Spanier \cite{Oldham-Spanier BOOK74},
Ross \cite{Ross BOOK75},
Samko, Kilbas \& Marichev   \cite{SKM BOOK93},
Podlubny \cite{Podlubny BOOK99}, 
Hilfer \cite{Hilfer BOOK00},
West, Bologna \& Grigolini \cite{West-Bologna-Grigolini BOOK03},
Kilbas, Srivastava \& Trujillo \cite{Kilbas-Srivastava-Trujillo BOOK06},
 Diethelm \cite{Diethelm BOOK10},
Mainardi \cite{Mainardi BOOK10}.
We have restricted our list  of  books published by 2010 because the list up to nowadays is too long, say impossible: any reader interested to Fractional Calculus
 can consult the  modern  books that appear  more suitable to his research field.
On the other side the Mittag-Leffler function has a less older story having been introduced by the Swedish mathematician Mittag-Leffler at the beginning of the 20-th century.
\\
{Many applications of Fractional Calculus and related special functions (e.g. the Mittag-Leffler function) can be found in the theory of viscoelasticity because of hereditary phenomena with long memory, with a story dating back to the twentieth century. Many authors indeed studied the properties of viscoelastic materials taking profits of the methods given by Fractional Calculus. Among this authors we can cite Scott Blair in  UK, Gemant in  USA and UK, Gerasimov and Rabotnov in Soviet Union. Among modern authors we can cite
 Meshkov and Rossikhin in Soviet Union,
 Caputo and Mainardi in Italy and later Torvik and Bagley in USA.
In particular, in 1971 Caputo and Mainardi  proved the role (and the presence) of the Mittag-Leffler function when derivatives of fractional order are introduced into the constitutive equations of a linear viscoelastic body.}\\
For more details {on the functions of the Mittag-Leffler type} the reader is referred to the treatise by
 Gorenflo et al \cite{GKMR BOOK14}
published in 2014 in the first edition. The second revised and enlarged edition
is planned by this year 2020.
\\
 For a minor acquaintance with this function the reader is referred to any book on fractional calculus for the great relevance that the Mittag-Leffler function has in this field.
 \vsp
{The idea of the paper is to revisit the previous works of Davis and Gross, bringing them to the attention of the scientific community.
 In particular, we have translated  from the original  German into English, the parts that were relevant for us, that, to the best of our knowledge, had not yet been done.
  The translation  is due for historical reasons,  because we wanted to give credits to these authors. However,  we have embedded  their results in a more detailed description of Fractional Calculus and related special functions in the light of  works, also by other authors, which have followed one another over the years.\\
Another name will appear in next sections, i.e. F.M. de Oliveira Castro, who worked in contact with Gross. In de Oliveira Castro's original work, that we 
have translated from German and reported here, the Mittag-Leffler function arises from the fundamental Volterra integro-differential equation of the problem making use of the classical method of successive approximations. 
In the present survey we show  that the same results can be obtained straightly  making use of the Laplace Transform, and this constitutes a simplification compared to the original work.\\
Our work can be seen from different points of view, being interdisciplinary.
Besides  revisiting and reviewing past results, our work deals with the solution
 of  problems of interest for physics and engineering making use of mathematical models and methods  in the spirit of modern  applied mathematics.}  
\vsp
  This paper is organized as follows. 
 In Section 2, we recall  some of the basic properties of the Mittag-Leffler function that is an  entire functions in the complex plane. 
In particular, restricting out attention in the non-negative time domain, we recall
for this function its completely-monotonicity properties,   
the corresponding  Laplace transforms and the  asymptotic expressions  
for small and large times. 
In such a domain   we show some plots in order to get a visualization of the Mittag-Leffler function.
 \\
In Section 3,   we devote our attention to the contributions of the mathematician Harold T. Davis who in the 1920's and 1930's recognized the role of the Fractional Calculus  with respect to the Volterra integral equation.
Then, as far as we know, he was the first to recognize the Mittag-Leffler function
in a contribution by Kennet S. Cole on nerve conduction. 
\\
 In Section 4, we outline the work of Bernhard Gross 
 who recognized the role of the Mittag-Leffler function in after-effect processes  in electrical circuits in 1930's and 1940's. Then, he devoted his attention to the mathematical  structure of linear viscoelastic bodies taking profit of the electrical mechanical analogy.  
\\
   Finally, in Section 4, we  provide some concluding remarks in order 
   to point out the long-standing  story of the applications of the 
   Mittag-Leffler functions due to their  scarce popularity 
in earlier times     in applied sciences.
   

\section{The Mittag-Leffler  functions}  
\noindent
The Mittag-Leffler function is defined by the following power series,
convergent in the whole complex plane,
$$
E_\alpha (z) := \sum_{n=0}^\infty \frac{z^n}{\Gamma (\alpha n+1)}
\,,\; \alpha > 0\,, \; z\in \CC\,. \eqno(2.1)$$
We recognize that it is an entire function   providing  a simple generalization of the
  exponential function to which it reduces for $\alpha=1$. 
 We also note that for the convergence of the power series in (2.1)
 the parameter $\alpha$ may be complex provided that $\Re (\alpha) >0$.
 Furthermore the function turns out to be entire of order 
 $1/\Re (\alpha)$.
The most interesting properties of the Mittag-Leffler function
are associated with its asymptotic expansions
as $z \to \infty$ in various sectors of the complex plane.
For  detailed asymptotic analysis,  
which includes the smooth transition across the Stokes lines,  
 the interested reader is referred to Paris \cite{Paris 2002},
  Wong and Zhao \cite{Wong-Zhao 2002}.
\\
In this paper we limit ourselves to the Mittag-Leffler function of order $\alpha \in (0,1)$ on the negative real semi-axis where is known to be completely monotone (CM) due a classical  result of 1948  by  Pollard 
\cite{Pollard 1948},
see also Feller
\cite{Feller BOOK71}.
\\
Let us recall that a function $ \phi(t)$  with $t\in\RR^+$ is    called a  completely
monotone (CM) function if it is non-negative, of class $C^{\infty}$,  and
$(-1)^n \phi^{(n)}(t)\ge 0$ for all $n \in \NN$. 
Then a function $ \psi(t)$  with $t\in\RR^+$ is    called a Bernstein function  
if it is non-negative, of class $C^{\infty}$, with a CM first derivative.
These functions play fundamental roles in linear hereditary mechanics to represent relaxation and creep processes, see e.g. Mainardi's book
\cite{Mainardi BOOK10}.
For mathematical details we refer the interested reader  to the survey paper by Miller and  Samko \cite{Miller-Samko 2001}
and to the more recent and exhaustive book by Schilling et al. 
\cite{Schilling BOOK12}.  
\\   
In particular we are interested for $t\ge0$ to the function 
$$ e_\alpha(t) := E_\alpha(-t^\alpha)
 =    
\sum_{n=0}^\infty (-1)^n \frac{t^{\alpha n}}{\Gamma (\alpha n+1)}
\,,\eqno (2.2)$$
with $0<\alpha \le 1$
that provides the solution to the fractional relaxation equation, 
see the 1997 survey  paper by Gorenflo and Mainardi (1997), 
see also Mainardi and  Gorenflo (2007), 
Mainardi \cite{Mainardi BOOK10}.  
\\
For readers'  convenience let us briefly outline the topic 
concerning the   generalization  via fractional calculus 
 of the first-order differential equation governing the
 phenomenon of (exponential) relaxation.  
Recalling   (in non-dimensional units)
 the initial value problem 
 $$\frac{du}{dt} = -u(t) \,, \; t\ge 0\,, \; \hbox{with}\; u(0^+)= 1\,\eqno(2.3)$$
 whose solution is
 $$u(t) = \exp (-t)\,,\eqno (2.4)$$
  the following  two alternatives 
  with $\alpha \in (0,1)$ are offered in the literature:
  $$\frac{du}{dt } = - D_t^{1-\alpha}\,u(t),\; t\ge 0, \;  u(0^+)= 1,\eqno(2.5a)$$
  $$ _*D_t^\alpha \,u(t) = -u(t), \; t\ge 0, \;  u(0^+)= 1,.
  \eqno(2.5b)$$
where  
 $D_t^{1-\alpha}$ and $\,  _*D_t^\alpha$ denote the fractional derivative of order $1-\alpha$ in 
 the Riemann-Liouville sense and the fractional derivative of order $\alpha$ in the Caputo sense, respectively.
 \\
 For a generic order $\mu\in (0,1)$ and for a sufficiently well-behaved function $f(t)$ with $t\in \RR^+$ 
 the above derivatives are defined as follows, see e.g. Gorenflo and Mainardi(1997), Podlubny (1999), 
$$ D_t^\mu  \,f(t) =
  {\ds  \rec{\Gamma(1-\mu )}}\, 
{\ds \frac{d}{dt}}\left[  
  \int_0^t
    \! \frac{f(\tau)}{ (t-\tau )^{\mu }}\,d\tau\right]\,, \eqno(2.6a)  $$ 
$$  _*D_t^\mu  \,f(t) =
{\ds \rec{\Gamma(1-\mu )}}\int_0^t
    \! \frac{f^\prime (\tau)}{ (t-\tau )^{\mu }}\,d\tau\,. \eqno(2.6b)  $$   
  Between the two derivatives we have the  relationship 
   $$
   \begin{array}{ll}
	 {\ds _*D_t^\mu  \,f(t)} &= 
  {\ds D_t^\mu \,f(t) -  f(0^+)\, \frac{t^{-\mu}} {\Gamma(1-\mu )}}\\
  &  =
 D_t^\mu  \,\left[ f(t) - f(0^+) \right] \, .
\end{array} 
 \eqno(2.7)
   $$
Both derivatives in the limit $\mu \to 1^-$ reduce to the standard first derivative but for 
$\mu \to 0^+$ we have
 $$
  D_t^\mu f(t) \to f(t)\,,\;
    _*D_t^\mu f(t) \to f(t) - f(0^+)\,, \
   \eqno(2.8)$$  
 \\
 In analogy to the standard problem (2.3), we  solve the problems (2.5a) and (2.5b)
  with the Laplace transform technique,  using  the rules pertinent to the corresponding fractional derivatives.
 The problems (a) and (b) are equivalent since  the Laplace transform of the solution
 in both cases comes out as
 $$  \widetilde u(s) = \frac{s^{\alpha-1}}{s^\alpha +1}\,, \eqno(2.9) $$
that yields our function
$$ u(t) = e_\alpha(t) := E_\alpha(-t^\alpha)\,.\eqno(2.10)$$
The Laplace transform pair
$$       e_\alpha(t) \,\div \, \frac{s^{\alpha-1}}{s^\alpha +1}\,, \quad \alpha >0 \,, \eqno(2.11)$$
 can be proved   by transforming term by term the power series representation  of $e_\alpha(t)$
 in the R.H.S of  (2.2).
 Furthermore,    by
 anti-transforming the R.H.S  of (2.11)  by using the complex Bromwich formula, and taking into account  for $0<\alpha<1$ the contribution from  branch cut on the negative real semi-axis (the denominator $s^\alpha +1$  does nowhere vanish in the cut plane $-\pi \le \hbox{arg} \, s \le \pi $),  we get, 
 see also Gorenflo and Mainardi (1997),
$$ e_\alpha(t) = \int_0^\infty \!\!  \e^{-rt} K_\alpha(r)\, dr\,, \eqno(2.12) $$
where     
$$     
\begin{array}{ll}
 K_\alpha(r) &= 
 {\ds \mp \,\rec{\pi}\,    {\rm Im}\,
  \left\{ \left. \frac{s^{\alpha -1}} {s^\alpha +1}\right\vert_{{\ds s=r\, \e^{\pm i\pi}}} \right\}}
  \\
 &= {\ds \rec{\pi}\,
   \frac{ r^{\alpha -1}\, \sin \,(\alpha \pi)}
    {r^{2\alpha} + 2\, r^{\alpha} \, \cos \, (\alpha \pi) +1}} \ge 0\,.
    \end{array}
      \eqno(2.13)$$
    Since   $K_\alpha(r)$ is non-negative for all $r$ in the integral, the above formula proves 
    that $e_\alpha(t)$ is  CM function in view of the Bernstein theorem. This theorem provides a necessary and sufficient condition for a CM function as a real Laplace transform of a non-negative measure. 
\\
However, the CM property of $e_\alpha(t)$  can also be seen  as a consequence  of the result by  Pollard
because the transformation $x=t^\alpha$ is a Bernstein function for $\alpha \in (0,1)$.
In fact it is known that a CM function can be obtained by composing a CM with a Bernstein function based on the following theorem:
{\it Let $\phi(t)$ be a CM function and let $\psi(t)$ be a Bernstein function, then $\phi[\psi(t)]$  is a CM function.} 
 \\
  As a matter of fact,  $K_\alpha(r)$ provides an interesting  spectral representation of  $e_\alpha(t)$ in frequencies. With the change of variable $\tau=1/r$ we get the corresponding spectral representation in relaxation times, namely
$$
\begin{array}{ll}
 e_\alpha(t) &= \int_0^\infty \e^{-t/\tau} H_\alpha(\tau)\, d\tau\,, \\
H_\alpha(\tau) &= \tau^{-2} \, K_\alpha(1/\tau)\,, 
\end{array}
\eqno(2.14) 
$$
that can be interpreted as a continuous distributions of elementary (i.e. exponential) relaxation processes.
As a consequence    we get the identity between  the two spectral distributions, that is 
$$ H_\alpha(\tau)
= \rec{\pi}\,
   \frac{ \tau^{\alpha -1}\, \sin \,(\alpha \pi)}
    {\tau^{2\alpha} + 2\, \tau^{\alpha} \, \cos \, (\alpha \pi) +1} 
\,, \eqno(2.15)$$
a surprising fact pointed out in Linear Viscoelasticity by 
Mainardi in his 2010 book \cite{Mainardi BOOK10}.
This kind of universal/scaling  property seems a peculiar one for our Mittag-Leffler function $e_\alpha(t)$.
In Fig 1 we show $K_\alpha(r)$ for some values of the parameter $\alpha$.
 Of course for $\alpha=1$ the Mittag-Leffler function reduces to the exponential function $\exp(-t)$
 and the  corresponding spectral distribution is the    Dirac delta generalized function centred at $r=1$, namely $\delta(r-1)$.
\\
 \begin{center}
\includegraphics[width=7.5cm]{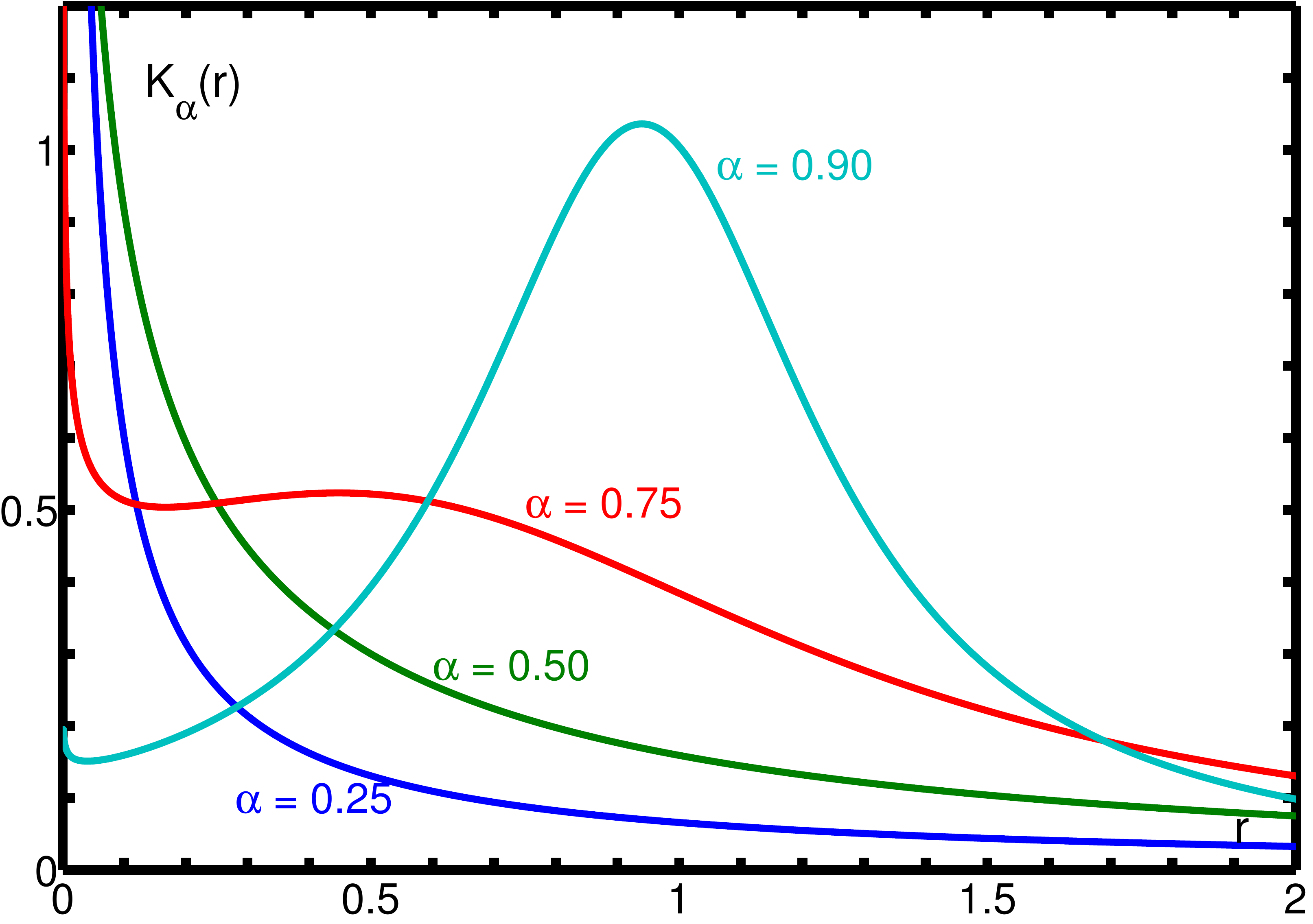}
\end{center}
{{\bf Fig.1} Plots of the spectral function  $K_\alpha(r) $ for $\alpha=0.25, 0.50, 0.75, 0.90$ in the  range  $0\le r \le 2$.}
\\
In Fig 2  we show some plots of $e_\alpha(t)$  for some values of the parameter $\alpha$.
It is worth to note the different rates of decay of $e_\alpha(t) $   for small and large times.
In fact the decay is very fast as $t\to 0^+$ and very slow as $t \to +\infty$.
\\ 
\begin{center}.
\includegraphics[width=7.5cm]{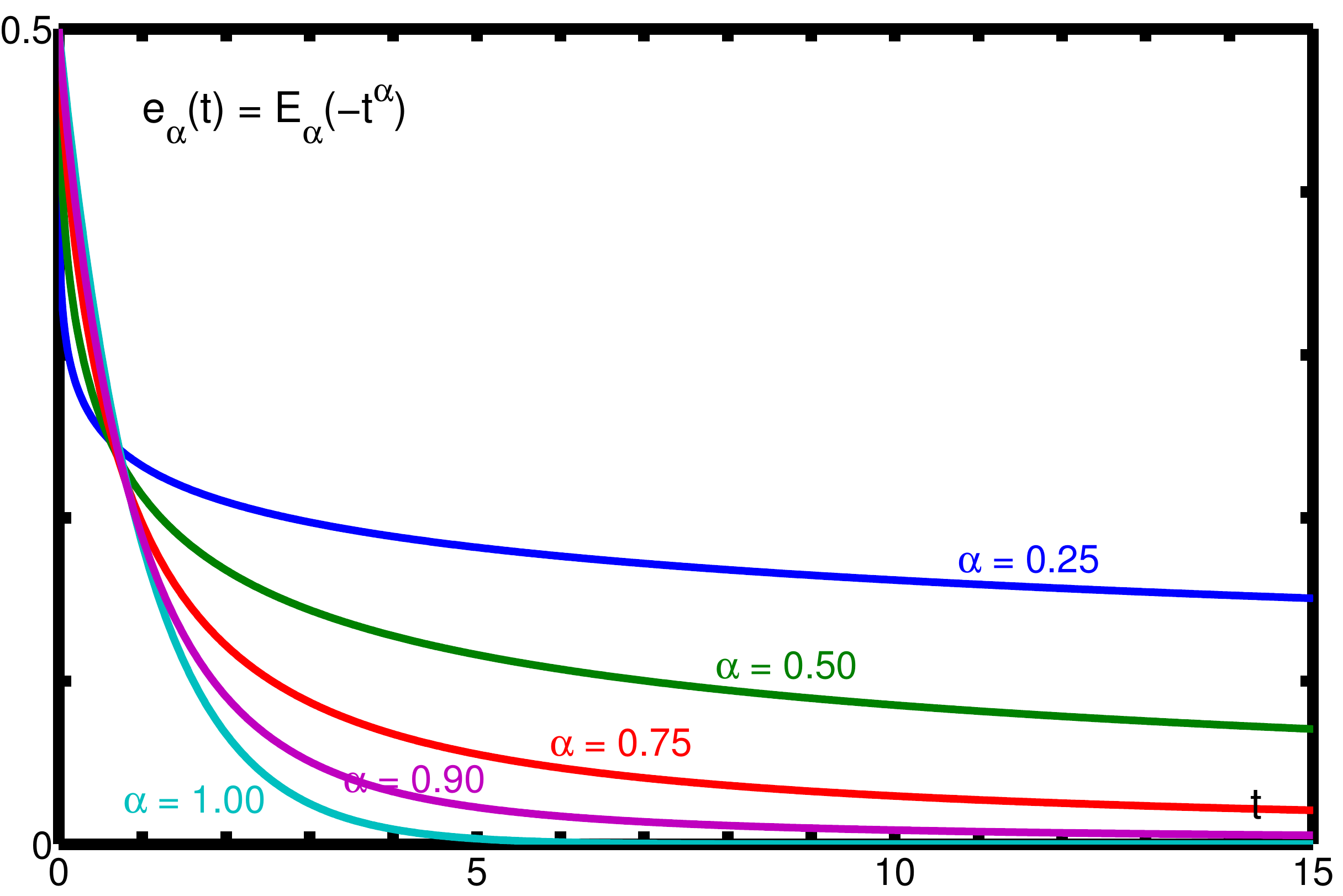} 
\end{center}
{{\bf Fig.2} Plots of the Mittag-Leffler function $e_\alpha(t)$  for $\alpha=0.25, 0.50, 0.75, 0.90, 1.$ in the  range $0\le t\le 15$ }
\\
\subsection{The two common asymptotic approximations}
  It is common to point out that     the function $e_\alpha(t)$  matches for $t\to 0^+$ with a stretched exponential    with an infinite negative derivative,   whereas as $t\to \infty$ with a negative power law.
The short time  approximation is  derived from the convergent power series representation (2.2). In  fact for $t\ge0$,
$$ e_\alpha(t) = 1 - \frac{t^\alpha}{\Gamma(1+\alpha)} + \dots 
\sim  \exp{\ds \left[- \frac{t^\alpha}{\Gamma(1+\alpha )}\right]}\,.
. \eqno(2.16)$$
The long time approximation is derived from the asymptotic power series representation of $e_\alpha(t)$ that turns out to be, see Erd\'elyi (1955),
 $$ e_\alpha(t) \sim 
 \sum_{n=1}^\infty  (-1)^{n-1} \,\frac{ t^{-\alpha n}}{\Gamma(1- \alpha n)}\,, 
 \; t\to \infty\,,\eqno(2.17) $$
 so that, at the first order, 
 $$e_\alpha(t) \sim  {\ds \frac{t^{-\alpha}}{\Gamma(1-\alpha )} }\,, \; t\to \infty\,. \eqno(2.18) $$
{The asymptotic representations (2.16)  for $t\to0$ and 
(2.17)--(2.18)  for $t\to \infty$ can be obtained by applying the Tauberian theory to the Laplace transforms (2.11)
for $s \rightarrow \infty$ and $s \rightarrow 0$, respectively.}
\\
 As a consequence  the function $e_\alpha(t)$  interpolates 
 for intermediate time $t$ between the stretched exponential
and the negative power law.
The stretched exponential models   
  the  very fast decay  for small  time $t$, whereas the asymptotic   power law 
  is due to the very slow decay for large  time $t$.
In fact, we have the two commonly stated  asymptotic  representations: 
  $$
\e_\alpha  (t) \sim
\left\{ 
\begin{array}{ll}
 e_\alpha^0(t) := \exp{\ds \left[- \frac{t^\alpha}{\Gamma(1+\alpha )}\right]}\,, 
  &  t\to 0\,;   \\ \\
e_\alpha^\infty(t) := {\ds \frac{t^{-\alpha}}{\Gamma(1-\alpha)}},
 &  t\to \infty\,.
\end{array}
\right . \eqno(2.19)
   $$
   The stretched exponential   replaces  the rapidly decreasing expression   
    $1- {t^\alpha}/{\Gamma(1+\alpha )} $  from Eq. 2.16.
 Of course, {\it for sufficiently small and for sufficiently large values of $t$} 
  we have the inequality
   $$ e_\alpha^0(t) \le e_\alpha^\infty (t)\,,  \quad 0<\alpha<1\,. \eqno (2.20)$$
   \\
In  Figs  3-7,
for $\alpha = 0.25, 0.5, 0.75, 0.90, 0.99$,
 we compare  in   logarithmic scales  the function $e_\alpha(t)$ (continuous line) and its asymptotic representations, the stretched exponential $e^0_\alpha(t)$ valid for $ t\to 0$ (dashed line)  and the power law $e^\infty_\alpha(t)$  valid for $t\to \infty$ (dotted line).
We have chosen the time range  $10^{-5} \le t \le 10^{+5}$.

\begin{center}
\includegraphics[width=7.35cm]{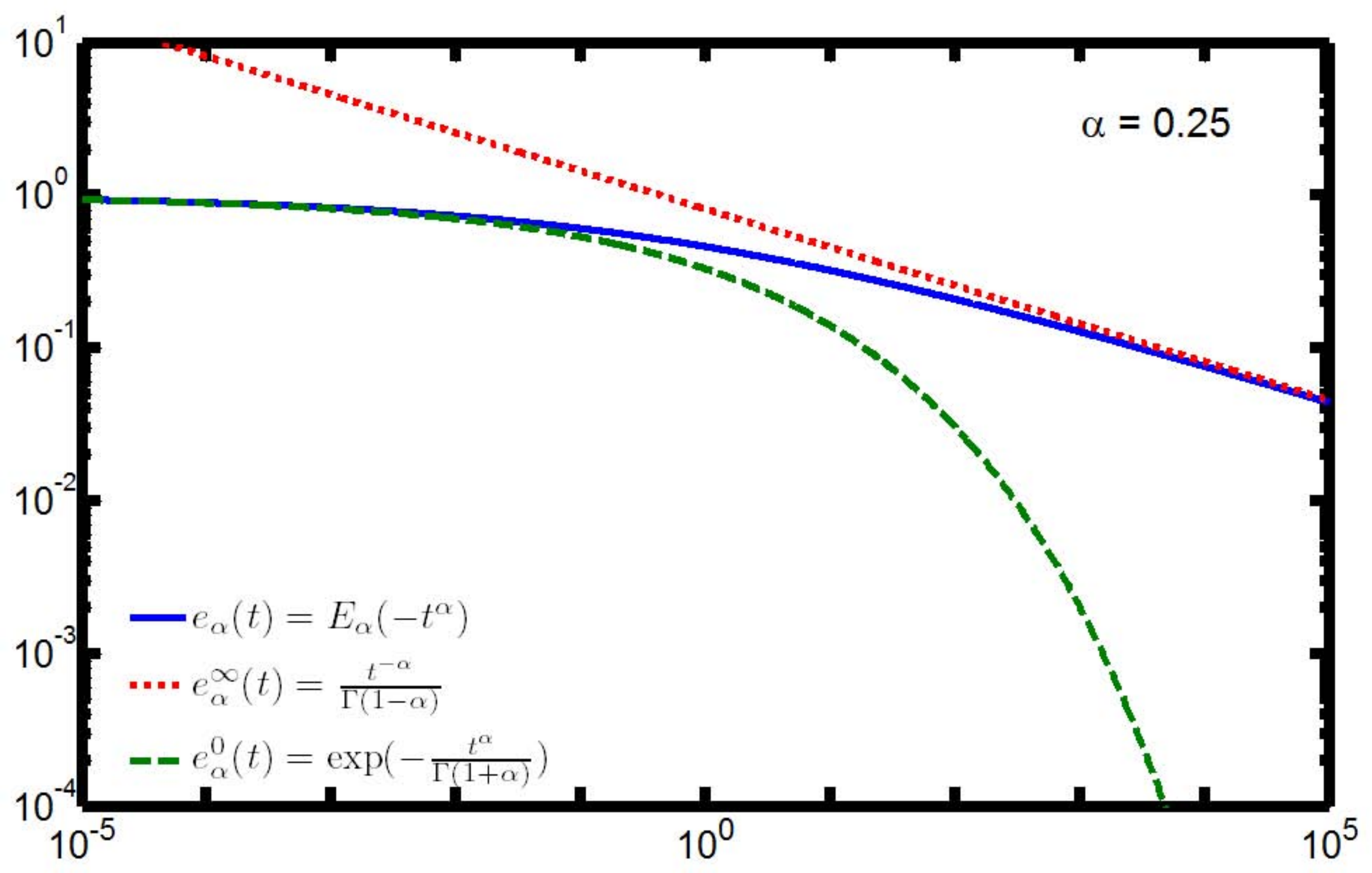}
\end{center}
\vskip-0.25truecm
{{\bf Fig.3} Approximations $e^0_\alpha(t)$ (dashed line) and 
$ e^\infty_\alpha(t)$ 
(dotted line) to   $e_\alpha(t)$ 
 in   $10^{-5} \le t \le 10^{+5}$;  $\alpha=0.25$.}
\smallskip
\begin{center}
\includegraphics[width=7.35cm]{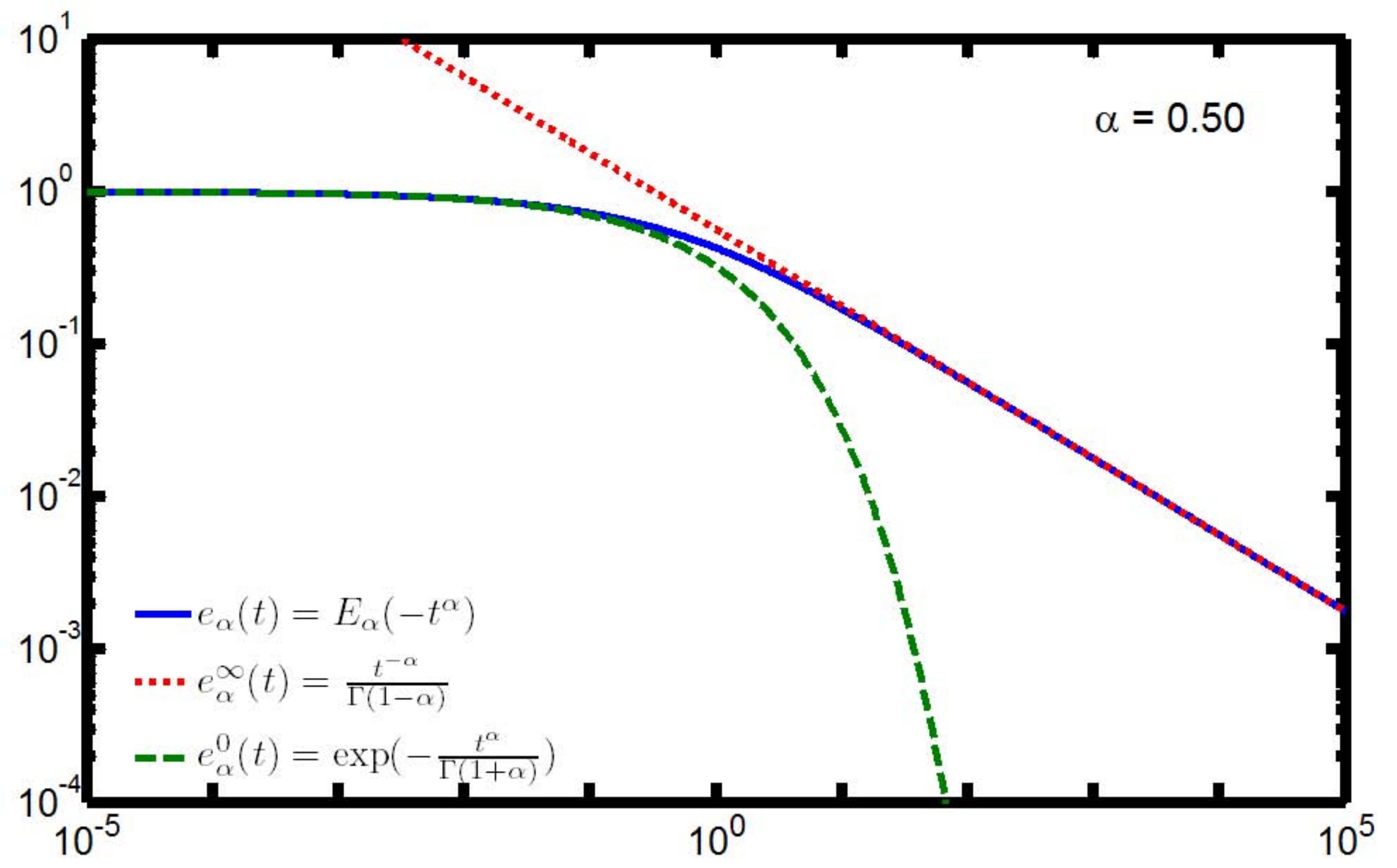}
\end{center}
\vskip-0.25truecm
{{\bf Fig.4} Approximations $e^0_\alpha(t)$ (dashed line) and 
$ e^\infty_\alpha(t)$  (dotted line) to  $e_\alpha(t)$     
 in  $10^{-5} \le t \le 10^{+5}$; $\alpha=0.50$.}
\smallskip
\begin{center}
\includegraphics[width=7.35cm]{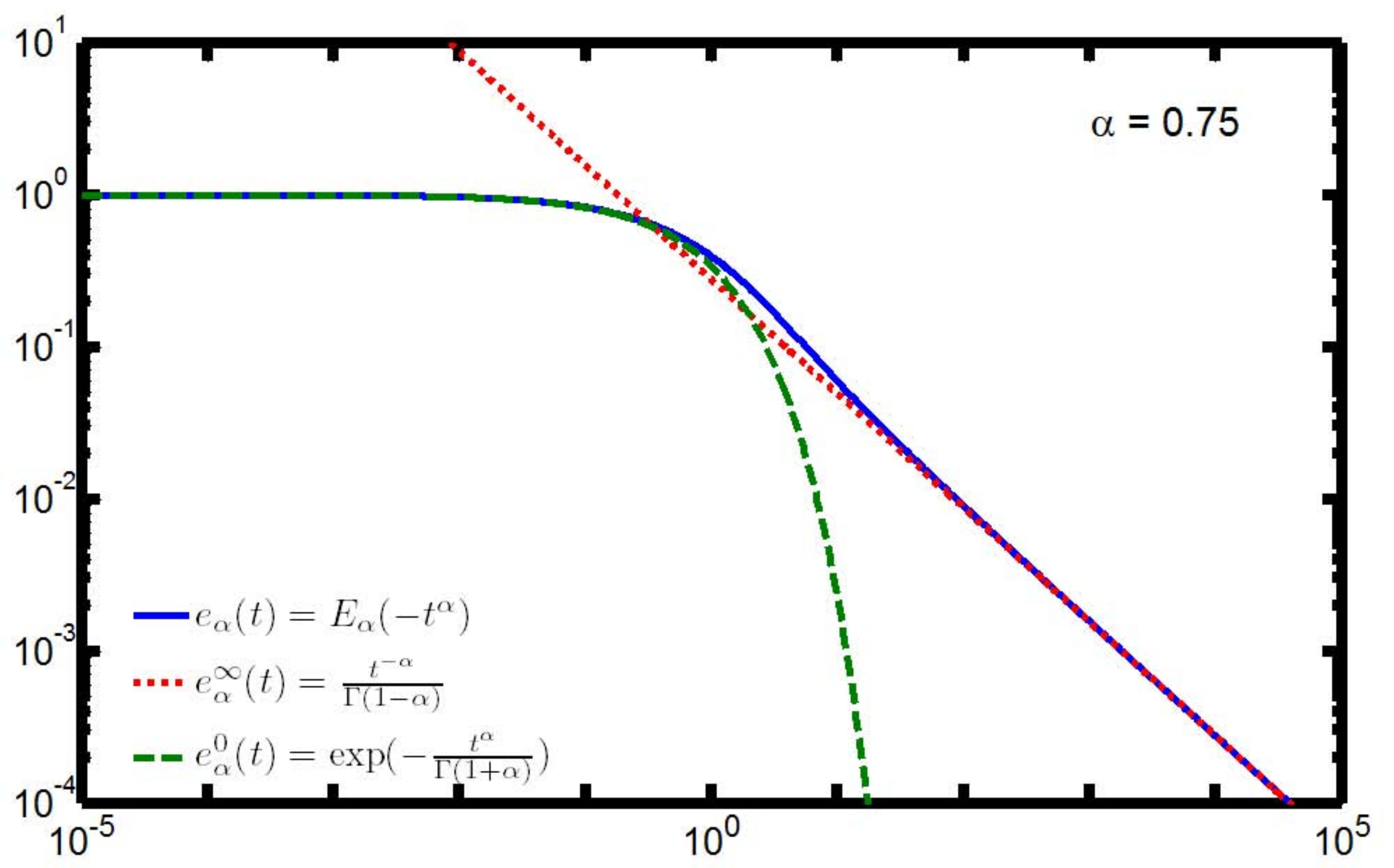}
\end{center}
\vskip-0.25truecm
{{\bf Fig.5} Approximations $e^0_\alpha(t)$ (dashed line) and 
$ e^\infty_\alpha(t)$ 
(dotted line) to   $e_\alpha(t)$    
 in   $10^{-5} \le t \le 10^{+5}$; $\alpha=0.75$.}
\smallskip
\begin{center}
\includegraphics[width=7.35cm]{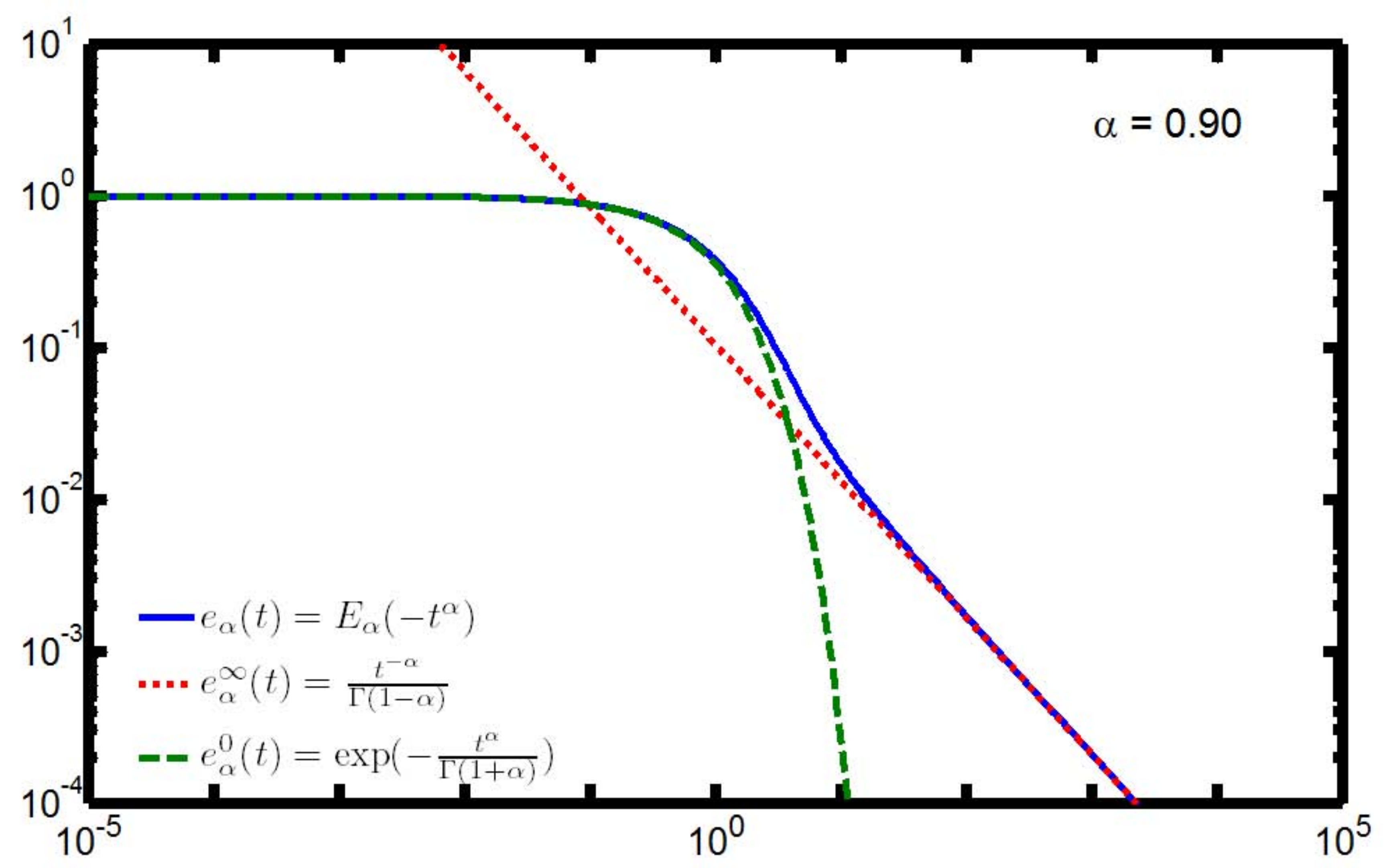}
\end{center}
\vskip-0.25truecm
{{\bf Fig.6} Approximations $e^0_\alpha(t)$ (dashed line) 
and $ e^\infty_\alpha(t)$ 
(dotted line) to   $e_\alpha(t)$    
 in   $10^{-5} \le t \le 10^{+5}$; $\alpha=0.90$.}
\smallskip
\begin{center}
\includegraphics[width=7.35cm]{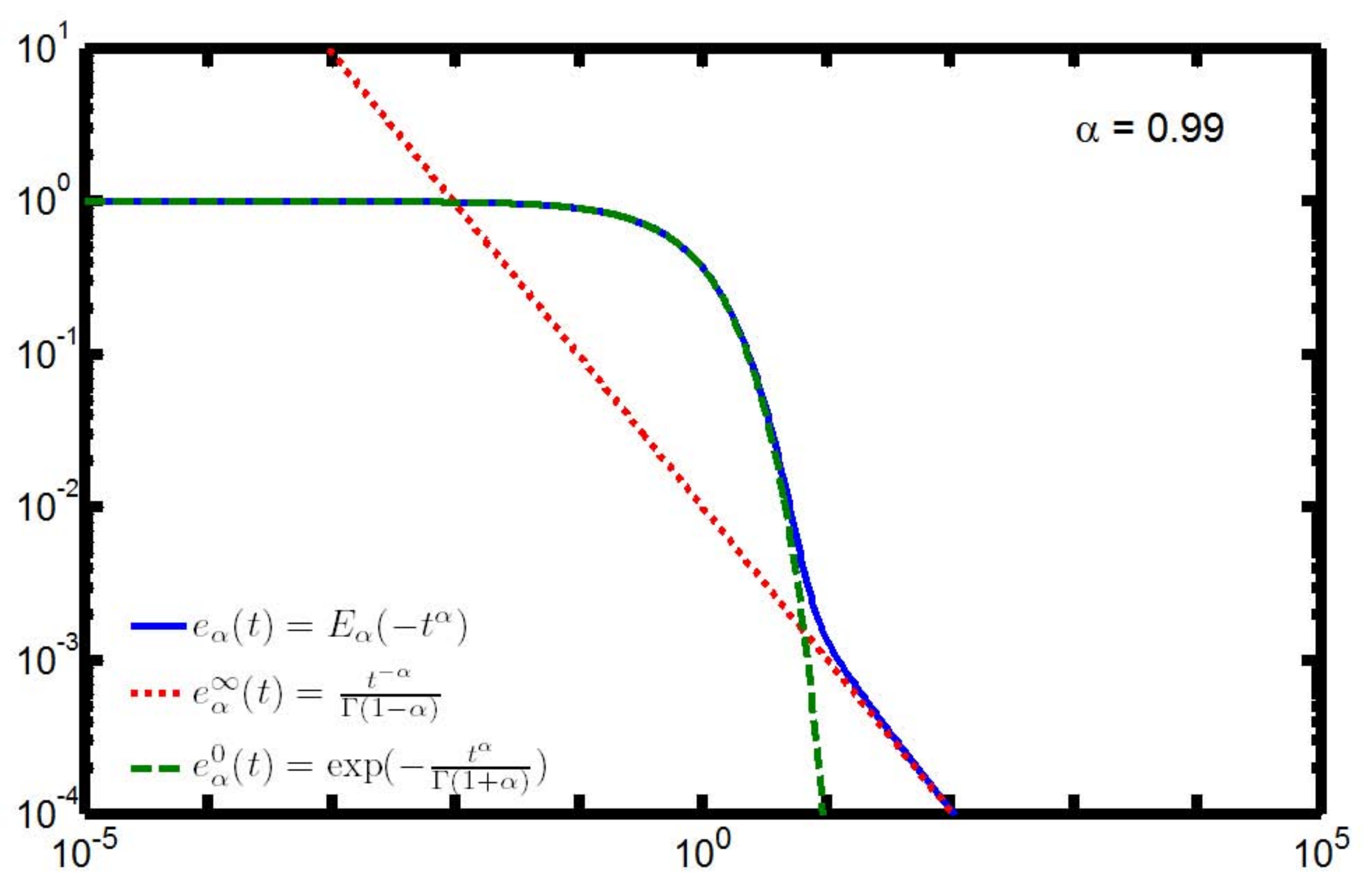}
\end{center}
\vskip-0.25truecm
{{\bf Fig.7} Approximations $e^0_\alpha(t)$ (dashed line) 
and $ e^\infty_\alpha(t)$ 
(dotted line) to   $e_\alpha(t)$    
 in $10^{-5} \le t \le 10^{+5}$; $\alpha=0.99$.}

In the 2014  paper by Mainardi 
\cite{Mainardi 2014}
the author, based on numerical computations,
 has conjectured the relation for $t>0$
 $$    f_\alpha(t)  \le  e_\alpha(t) \le g_\alpha(t)
 \,,
 \eqno(2.21)$$ 
 where
   $$ 
\begin{array} {ll}   
   f_\alpha(t) &:= {\ds \frac{1}{1+ \frac{t^\alpha}{\Gamma(1+\alpha)}}}, \\
    g_\alpha(t) &:= {\ds \frac{1}{1+ t^\alpha \Gamma(1-\alpha)}}\,.
    \end{array}
    \eqno(2.22)$$
 This conjecture was proved by Simon with rigorous arguments based on the probability theory, see \cite{Simon 2014}.   
\section{The contributions  of Harold T. Davis}
In his 1936 book \cite{Davis BOOK36}
the mathematician Harold T. Davis devoted 3 sections (7, 8,9)  in Ch II, pp 64--76
to fractional operators, including fractional integrals and derivatives (in the sense of Riemann-Liouville). Then he devoted
one section (7) in Ch VI, pp 280-294  on special  applications of the 
Fractional Calculus mostly involving the Mittag-Leffler function.
\\
For the fundamentals of the fractional calculus  we recall the previous notes by Davis published in 1924
 \cite{Davis 1924}
 and in 1927 \cite{Davis 1927}.
\\
For the applications of the Mittag-Leffler  we outline the 1930 fundamental paper by Hille and Tamarkin \cite{Hille-Tamarkin 1930}
where the authors solved the Abel integral equations of the second kind
by using the Mittag-Leffler function.
But this mathematical paper presumably was not observed by Davis 
who outlined in Ch VI different  applications of the Mittag-Leffler functions  in some examples.
Hereafter  we describe the example 4 where Davis reports on two notes by Kennet Cole on nerve conduction published 1933 in the proceedings
of the First Symposium on Quantitative Biology (Cold Spring Harbor),
see \cite{Cole 1933}.
This report gave to Davis the occasion to interpret the solution of Cole
in terms of the Mittag-Leffler function of which Cole (being 
a physicist interested to biology) was not aware.
    In our opinion, this result may be considered the first appearance of the 
    Mittag-Leffler functions in applied sciences.
\subsection{The Mittag-Leffler function in nerve \\ conduction}
 A living nerve can be stimulated by passing a direct current
through a short portion of it between two electrodes, provided the
potential difference exceeds a certain critical value known as the
rheobase. As the duration of the potential applied across the 
electrodes decreases, it is found that the intensity necessary for stimulation
increases rapidly in a hyperbolic manner. The following analysis
is designed to explain this phenomenon.
\\
An idealized nerve fiber consists of a cylindrical core of electrolyte
covered with a thin sheath or membrane. It is assumed that a
local threshold change of the normal potential difference across the
membrane will stimulate the fiber and cause an impulse to be propagated.
\\
The problem is then to express analytically the strength of
stimulus which, when applied to the nerve bundle as a whole, will
change the potential difference across the membrane of an individual
fiber by a threshold amount in a given time
 \\
To begin with, experimental evidence points to the conclusion
that the electrical behavior of the nerve fiber may be simulated by
the type of circuit illustrated in Fig. 8, where $r$ and $R$ are constant
resistances and the element $P$, called the polarization element,
has an impedance  defined by the following equation :
$$
K p^{-\alpha}  \to I_P(t) = e_p(t)\ \quad
0< \alpha <1. 
\eqno(3.1)
$$
$I_r (t)$ and $I_P (t)$ are the instantaneous currents in $r$ and 
$P$.
Then 
$e_p (t)$ is  the potential
across the  element $P$  and $p$ denotes  operator $d/dt$. 
The positive constant $K$ is
determined experimentally. 
\\
No combination of electrical circuits with
ordinary resistances and capacities is known to lead to an impedance
of the form postulated, but experimental evidence appears to indicate
that such an impedance is essential to the description of the
curious electrical behavior of biological materials in general and of
nerve fibers in particular.
\begin{center}
\includegraphics[width=6.00cm]{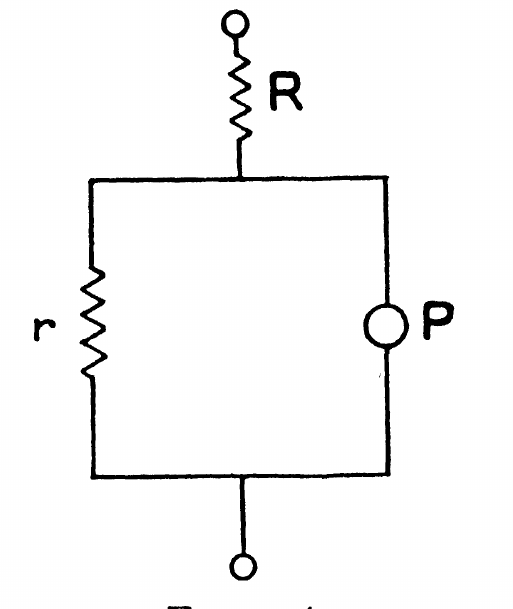}
\end{center}
\vskip-0.25truecm
{{\bf Fig. 8} The Cole circuit {useful to simulate the electrical behavior of a nerve fiber consisting of a cylindrical core of electrolyte covered with a thin sheath or membrane.}} 
\\
Referring now to the figure 8, we compute the relationship between
$I_p (t)$ and $e_p (t)$, when a constant electromotive force (e. m. f.)
$E$, is applied across the electrodes.
Applying the Kirkoff laws for the circuits and applying the Laplace transform 
the solution is easily found to be
$$  e_P(t) = \frac{ER}{R+r}
\left [ 1 - E_\alpha(-\lambda t^\alpha) \right],  
\eqno(3.2)
$$
where
$$ \lambda =  K\, (R+r)/Rr\,.\eqno(3.3)$$
Using non dimensional units for the current and the time, we 
have the following visualization of Eqs (3.2)-(3.3) at variance of
selected values  the parameter $\alpha\in (0,1]$.
We note that for $\alpha=1$ the Mittag-Leffler reduces to the exponential
 function with a great difference with the cases $\alpha \min (0.1)$
 because of the strong variability for small times ad low variabilty for large times as expected from the corresponding asymptotic expressions stated in Eq (2.19).     
\begin{center}
\includegraphics[width=8.5cm]{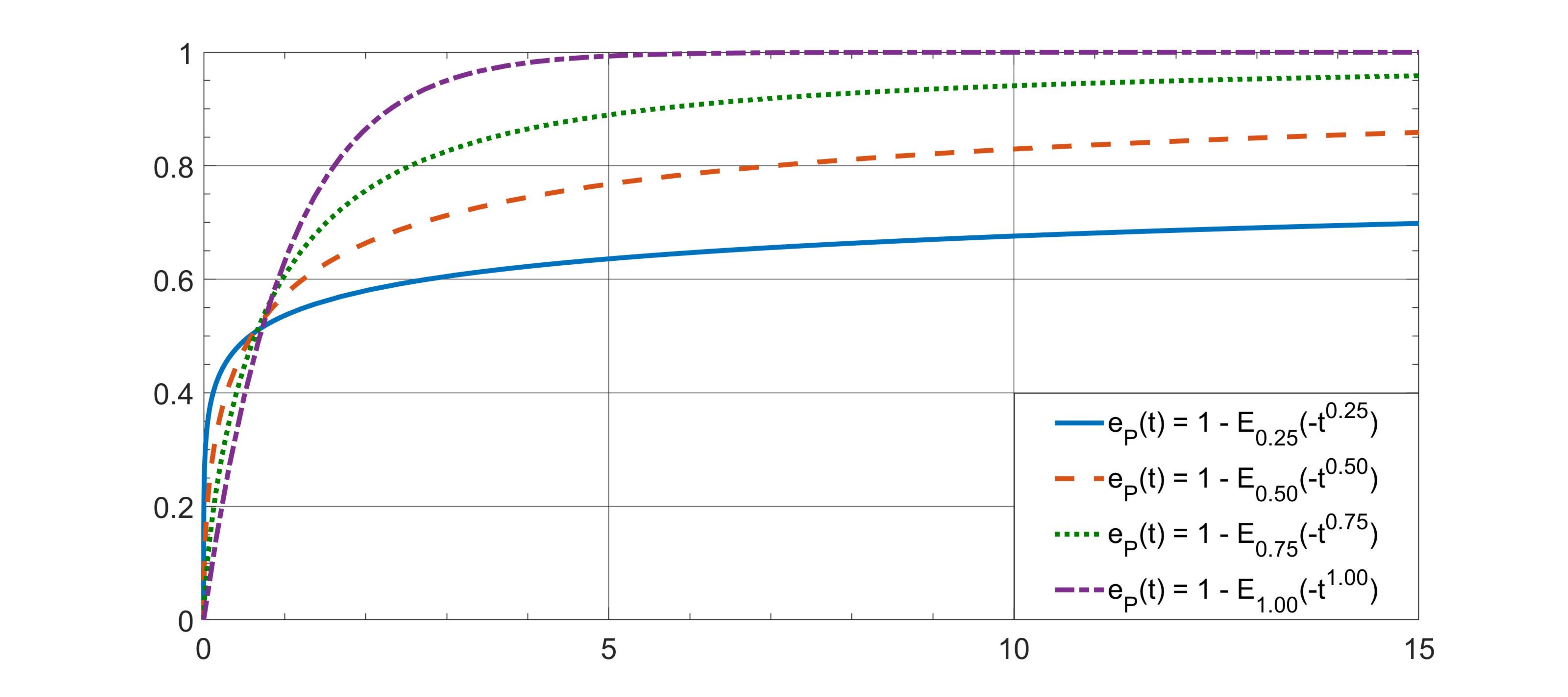}
\end{center}
\vskip-0.25truecm
{{\bf Fig.9} The potential across the element $P$ in the Cole circuit for nerve conduction, versus time,  for selected values of $\alpha= 0.25,~0.50,~0.75,~1.00.$} 
\smallskip
\\
We then outline that Kennet Cole with his younger brother Robert
published two seminal papers on dielectric relaxation,
\cite{Cole-Cole 1941,Cole-Cole 1942} illustrated in the 
historical survey by Valerio, Machado and Kiryakova
\cite{Valerio-Machado-Kiryakova 2014} on some pioneers of the applications of fractional calculus. Indeed the laws for dielectric relaxation that Cole brothers introduced and known as {\it Cole-Cole models} can be considered expressed in terms of fractional operators and Mittag-Leffler function, even if the authors were (once again) not aware of these interpretations.
\\
For more details on dielectric laws interpreted by  fractional operators and special functions of the Mittag-Leffler type,    we refer the reader to the more recent survey by 
Garrappa, Mainardi and Maione
\cite{Garrappa-Mainardi-Maione 2016}

 \section{The contributions of Bernhard  Gross}
Bernhard Gross was  a  physicist graduated in  Stuttgart, Germany.
Since  he arrived in Rio de Janeiro in Brazil in 1933, 
for many years he was   associated with an
electrical engineering department. 
He became familiar with the
mathematics of dielectric phenomena and electrical network
theory - the counterparts of viscoelastic effects - a long time
before he had his attention directed to the latter. This will
perhaps excuse his tendency of bringing into the discussion of
viscoelastic theory concepts and methods which originated in
electrical theory.
\\
After producing several papers in electrical circuits, 
Gross gave the opportunity
in 1939  to F.M de Oliveira Castro to discuss with more details
the mathematical theory of the  voltage   
of a capacitor which has a non exponential after-effects
\cite{Castro 1939}.
Indeed  Gross in 1937 \cite{Gross 1937}  had noted
that the voltage curve of a capacitor which has after-effects, when it is discharged, is not exponential and depends on the charging time; after a temporary short circuit, residual charges appear again. The discussion of Gross 
needed of a rigorous mathematical approach
based on Volterra integral equations that leads to Mittag-Leffler functions.
\\
In the following we report briefly the approach by Castro but we use the most
convenient tool of the Laplace transforms.  

\subsection{Introduction of the paper by Castro}
 In analogy to the simple differential equation of the capacitor without after-effects, an integro-differential equation can be set up if the superposition principle is valid for the capacitor with after-effects.  Gross has solved this approximately and the solution reproduced the observed behavior correctly in principle.  The degree of approximation could not be overlooked; so the question remained open whether existing deviations between the measurement and the calculation are due to inadequacies of the calculation or to the non-fulfillment of its requirements.\\
The aim of the present work is a strict integration of the discharge equation.  Their implementation provides a satisfactory theory of capacitor discharge;  moreover, we believe that it is a contribution to the theory of the Volterra integral equation, which is also of interest from a purely mathematical point of view.  The prerequisite for the calculation is the selection of a suitable after-effect function.  We have chosen the Schweidler expression $\beta t^{-n}$ for this, which has been shown time and time again to best represent the course of the recharging current in a wide interval \cite{Tank,Gnann,Wikstrom,HippaufStein,Voglis}; possible deviations for short times, such as Voglis \cite{Voglis} has determined is of no importance here.
\\
The work is organized as follows: the discharge equation is first converted to the form of a Volterra integral equation by a simple transformation.  The general integral of this equation is written and then solved via Laplace Transform technique; the solution can be simplified by limiting it to short discharge times, because here the solving core leads to a known transcendent, the function of Mittag-Leffler \cite{MittagLeffler 1905}. Finally, the solution obtained here is compared with the formulas given by Gross \cite{Gross 1937};  in the area for which they were derived, they prove to be a suitable approximation.
\subsection{The basic equation of the problem}
The assignments of an imperfect capacitor may be isolated at the instant $t = 0$. The voltage curve $U(t)$ is then calculated according to Gross \cite{Gross 1937}
 as an integral of the equation:
$$
C\frac{dU}{dt} + \frac{U}{R} + \int_{0}^{t}{\frac{dU}{d\tau} \phi (t - \tau)d\tau } + i_0 (t) = 0
\eqno(4.1)
\label{eq:integrodiffeq}
$$
Here $C$ is the geometric capacitance of the arrangement, $R$ is the ohmic terminal resistance, $\phi$ is the after-effect function and $i_0$ is the recharging current that results from all voltage changes during charging.\\
{In next subsections we will see how it is possible to apply the Laplace Transform to Eq. 4.1, after turning it into a Volterra integral equation.}\\
We restrict ourselves to the two particularly important cases of discharging after a charge that took place during the time $t_0$ with a voltage $U_0$, and recharging which is preceded by a full charge with $U_0$, and then by a short circuit during the time $t_0$.\\
Then we can write:
$$
i_0 (t) = \delta U_0 \phi (t + t_0)
\eqno(4.2)
$$
$$
[U(t)]_0 = U(0)
\eqno(4.3)
\label{eq:initialcondition}
$$
For the discharge is:
$$ 
	\delta = 1 ~~~~~~ \mathrm{and} ~~~~~~ U(0) = U_0
$$ 
and for recharge is:
$$ 
\delta = -1 ~~~~~~ \mathrm{and} ~~~~~~ U(0) = 0
$$ 
Eq. 4.3 represents the initial condition which has to be satisfied.\\
$\phi (t)$ is set, for the reasons given earlier, as:
$$
\phi (t) = \beta t^{-n}, ~~ (0 \leq n \leq 1, ~ \beta > 0)
\eqno(4.4)
$$

\subsection{Transformation of the basic equation into a Volterra integral equation of the second kind}
Eq. 4.1 is an integro-differential equation for $U(t)$, which can be immediately brought to the well-known form of a Volterra integral equation. It is introduced the unknown:
$$
\psi (t) = \frac{dU}{dt}
\eqno(4.5)
$$
obtaining
$$
\psi (t) + \int_{0}^{t}{\psi (\tau) K(t - \tau) d\tau} = f(t).
\eqno(4.6)
\label{eq:integrodifferentialequation}
$$
The kernel of this equation is:
$$
K(t - \tau) = \lambda [1 + k(t - \tau)^{p-1}].
\eqno(4.7)
$$
The right term is:
$$
f(t) = -\bigg[\lambda U(0) + \frac{i_0 (t)}{C}\bigg]
\eqno(4.8)
$$
The function $f$ is restricted in the interval $0 \leq t \leq a$ ($a$ finite), when $0 < c \leq t_0$. In the above equations this abbreviations are used:
$$
k = \beta R, ~~~~~~ p = 1 - n, ~~~~~~ \lambda = 1/RC
\eqno(4.9)
$$
In the most general case, the discussion of the form of the solution and its application for performing numerical calculations is cumbersome; let us so first deal with two particularly important cases.
\subsection{Case a: $\mathbf{R = \infty}$; very short \\ charging or discharging times}
If the charging and discharging time is very short, experience has shown that the anomalous current component initially outweighs the ohmic component over a considerable time interval. In this interval $U/R$ can be neglected compared to $i_0 (t)$. Formally, this is done by setting $R = \infty$.
\\
Here we propose a Laplace Transform based method to solve Eq. 4.6. The sign $\div$ is used to indicate the juxtaposition among a function and its Laplace Transform, so that we have:
$$ 
	\psi (t) \div \mathcal{L}[\psi (t)] = \tilde{\psi}(s)
$$ 
$$ 
f(t) \div \mathcal{L}[f(t)] = \tilde{f}(s)
$$ 
$$ 
K (t) \div \mathcal{L}[K (t)] = \tilde{K}(s)
$$ 
It is straightforward to obtain, from Eq. 4.6:
$$
\tilde{\psi}(s) + \tilde{\psi}(s)\times \tilde{K}(s) = \tilde{f}(s)
\eqno(4.10)
$$
Hence:
$$ 
	\tilde{\psi}(s)[1 + \tilde{K}(s)] = \tilde{f}(s)
$$ 
$$ 
\mathcal{L}\bigg[ \frac{dU}{dt}\bigg] = \tilde{\psi}(s) = \frac{\tilde{f}(s)}{1 + \tilde{K}(s)}
$$ 
$$ 
	s \tilde{U}(s) - U(0) = \frac{1}{1 + \tilde{K}(s)}\tilde{f}(s)
$$ 
$$
\tilde{U}(s) = \frac{U(0)}{s} + \frac{1}{s(1 + \tilde{K}(s))}\tilde{f}(s)
\eqno(4.11)
\label{eq:LTpotential}
$$
Note that:
$$
\tilde{K}(s) = \frac{\lambda}{s} + \frac{\beta}{C}\frac{\Gamma (p)}{s^p}
\eqno(4.12)
$$
The latter expression is simplified as $R \rightarrow \infty$ ($\lambda \rightarrow 0$), and we substitute it inside Eq. 4.11. In particular:
$$
\frac{1}{s(1 + \tilde{K}(s))} 	= \frac{s^{p-1}}{s^p + \frac{\beta}{C}\Gamma (p)} = \tilde{e}(s)
\eqno(4.13)
\label{eq:LTMittagLeffler}
$$
The inverse Laplace Transform of Eq. 4.13 gives us the Mittag-Leffler function:
$$
\tilde{e}(s) \div \mathcal{L}^{-1}[\tilde{e}(s)] = E_p \bigg[-\frac{\beta}{C}\Gamma (p) t^p \bigg]
\eqno(4.14)
$$
\\
Thus, the solution is:
$$
U(t) = U(0) + \int_{0}^{t}{E_p \bigg[-\frac{\beta}{C}\Gamma (p) (t-s)^p  \bigg] f(s) ds}
\eqno(4.15)
\label{eq:solution}
$$
with:
$$
f(s) = -\frac{i_0}{C}
\eqno(4.16)
$$
Now the second special case will be dealt, which leads formally to the same expression as the above. 

\subsection{Case b: [U - U (0)]/R = 0;  shortly \\
after the opening}
In Eq. 4.1, instead of $U/R$, we write the expression $U/R - U(0)/R + U(0)/R$. $U(0)/R$ is constant and we can incorporate it in the right member; this then only changes the value of $f(s)$. If we limit ourself to such short times that $\frac{U - U(0)}{R}$ can be overlooked, we formally fall back to case $a$. The solution is again given by Eq. 4.15; the only difference is that now:
$$
f(s) = -\bigg(\frac{U(0)}{RC} + \frac{i_0 (s)}{C}\bigg)
\eqno(4.17)
\label{eq:formula35deoliveiracastro}
$$
\subsection{Calculation of the Mittag-Leffler \\ function}
For some special values $E_p$ leads to simple expressions. Indeed:
$$
E_1 (-x) = e^{-x},
\eqno(4.18)	
$$
$$
E_{0.5} (-x) = e^{x^2} [1 - \Phi (x)],
\eqno(4.19)
$$
where $\Phi (x)$ is the Gaussian error integral: 
$$ 
	\Phi (x) = \frac{2}{\sqrt{\pi}}\int_{0}^{x^2}{e^{-s^2}ds}
$$ 
Furthermore we have:
$$
~~~~~~~~~~~~~~~~	E_0 (-x) = \frac{1}{1 + x}, ~~~~~~~~~~ \mathrm{for} ~ |x| < 1
\eqno(4.20)
\label{eq:formula38deoliveiracastro}
$$
For $|x| > 1$ the function $E$ is undefined. However, the curve $\frac{1}{1 + x}$ still seems to represent a limit curve for the function $E_p$ if $p \rightarrow 0$; we do not conclude this on the basis of strict proof, but on the basis of numerical agreement, as we will shortly show.\\
From an experimental point of view small values of $p$, of the order of 0.1, are of interest. We have therefore calculated the function $E_{0.1}$ in the range of $x$ in question.\\
For $x \leq 1$ the calculation can be done with the following series:
$$
E_p (x) = 1 + \frac{x}{\Gamma (p + 1)} + \frac{x^2}{\Gamma (2p + 1)} + ... + \frac{x^h}{\Gamma (hp + 1)} + ...
\eqno(4.21)
\label{eq:formula31deoliveiracastro}
$$
For $x > 1$, instead, this is very tedious.  However, there are asymptotic formulas that are very already convenient for values of $x$ around 2 or more. The derivation should take place elsewhere.\\
We have, with $p = 1/m$ and $m$ being an even number
for  $ x \gg 1$:
$$
\begin{array}{ll}
\hskip-0.25truecm  &
E_p (-x)  =
{\ds \sum_{\nu = 1}^{n-1}{\frac{(-1)^{\nu + 1}}{\Gamma (1 - \nu /m)x^{\nu}}}}\\
&{\ds \bigg(1 - \frac{\nu /m}{x^m} + \frac{\nu /m (\nu /m + 1)}{x^{2m}} - ...\bigg), }
\end{array}
\eqno(4.22)
\label{eq:formula39deoliveiracastro}
$$
It is easy to see that the formula usually converges very quickly. In general, the term with $x^m$ is negligible compared to 1. Table 1 gives the values calculated in this way.\\
For the further execution of the calculation, which still requires an integration on $E$, it would be desirable to find a simpler even if only approximately valid representation of the function $E$.\\
As Eq. 4.20 is valid for $p = 0$, it is natural to try a  generalized approach for small values of $p$
$$
E_p (-x) = \frac{1}{1 + ax}
\eqno(4.23)
\label{eq:formula40deoliveiracastro}
$$
One can then determine $a$ so that the slope at the origin is correctly represented, and because of $\Gamma (1+p) = p\Gamma (p) $ it follows from Eq. 4.21 $a = 1/p\Gamma (p)$, and thus here:
$$
~~~~~~~~~~~~~~~~	E_p (-x) = \frac{1}{1 + x/p \Gamma (p)}, ~~~~~ p \ll 1
\eqno(4.24)
\label{eq:formula41deoliveiracastro}
$$
The values ​​calculated in this way can also be found in Table 1.
In particular, the calculations are performed for $p = 0.1$.\\
We can see that the match is very good, and it decreases with increasing $x$. Using the Eq. 4.22 the deviation can be estimated.\\
For $x \rightarrow \infty$ is given by first approximation given by $1/\Gamma (1 - p) x$, while Eq. 4.24 gives the expression $\Gamma (1 + p)/x$. The difference is meaningless for very small values of $p$.\\
\begin{center}
\begin{tabular}{ |c|c|c|c| } 
	\hline
	$\mathbf{x}$& $\mathbf{E_{0.1} (-x)}$ & $\mathbf{1/1 + 1.051x}$ & $\mathbf{E_{0.5} (-x)}$\\
	\hline 
	0.0 & 1.000 & 1.000 & 1.000\\ 
	\hline
	0.2 & 0.8259 & 0.8264 & 0.8090\\ 
	\hline
	0.4 & 0.7031 & 0.7040 & 0.6708\\ 
	\hline
	0.6 & 0.6118 & 0.6133 & 0.5678\\ 
	\hline
	1.0 & 0.4856 & 0.4876 & 0.4276\\ 
	\hline
	2.0 & 0.3200 & 0.3224 & 0.2655\\ 
	\hline
	4.0 & 0.1901 & 0.1922 & 0.1370\\ 
	\hline
	6.0 & 0.1353 & 0.1369 & 0.0940\\ 
	\hline
	8.0 & 0.1049 & 0.1063 & 0.0650\\ 
	\hline
	10.0 & 0.0857 & 0.0869 & 0.0564\\ 
	\hline
\end{tabular}
\end{center}
If the function $E$ can be replaced by Eq. 4.24, the calculation is further simplified.
\subsection{Explicit expression of the solution for short discharge times in the limit case of full charge}
We're considering now $p \ll 1$ and the case of discharge after a full charge. Then if
we limit ourselves to short discharge times, $U$ is given by Eqs. 4.15, 4.17 and 4.23:
$$
U(t) = U_0 - \frac{U_0}{RC}\int_{0}^{t}{\frac{ds}{1 + \frac{\beta}{pC}s^p}}
\eqno(4.25)
\label{eq:formula42deoliveiracastro}
$$
If $p$ can be represented as $1/m$, with $m$ being even, the integral in Eq. 4.29 can be evaluated.\\
With:
$$
A = \frac{\beta}{pC}
\eqno(4.26)
$$
the integral is written as:
$$
J = \frac{1}{p A^m}\int_{0}^{At^p}{\frac{u^{m-1} }{1 + u}du}.
\eqno(4.27)
$$
By dividing and integrating it finally follows:
$$ 
\begin{array}{ll}
\hskip -0.5truecm
	J = {\ds \frac{1}{p A^m}}
\hskip -0.2truecm	
	& {\ds \bigg[\frac{A^{m-1} t^{1-p}}{m-1} - \frac{A^{m-2} t^{1-2p}}{m-2} + ... - ..}. \\
	&{\ds  + At^{1 - (m-1)p} - \ln{(1 + At^p)}\bigg]}
\end{array}	
\eqno(4.28)
$$
This case is completely solved.
\subsection{Comparison between the rigorous \\ solution and Gross's expressions}
Gross \cite{Gross 1937} obtained an approximate solution to Eq. 4.1 on the assumption that $dU/d\tau$ changes very slowly compared to $\phi (t - \tau)$, and therefore can be taken out
of the integral.\\
Then, in our notation:
$$
U(t) = U_0 \exp{\frac{1}{U_0}\int_{0}^{t}{\frac{f(s)ds}{1 + \frac{\beta C}{p}s^p}}}.
\eqno(4.29)
\label{eq:formula49deoliveiracastro}
$$
In the two special cases and for $p \ll 1$, the strict solution could be written in the form:
$$
U(t) = U_0 + \int_{0}^{t}{\frac{f(s) ds}{1 + \frac{\beta C}{p}(t-s)^p}}.
\eqno(4.30)
\label{eq:formula50deoliveiracastro}
$$
This already takes the form of Eq. 4.19. A first condition for the approximation to be valid is therefore given by the requirement $p \ll 1$, or $n \approx 1$. But that is exactly the assumption.\\
We now want to compare the solutions more closely.
\\
a) In the case of very long loading times and discharge times so short that we can truncate Eq. 4.29 at the linear term, the two equations are identical. This is because here, according to Eq. 4.17, $f(s) = -U_0 /RC = \mathrm{const}$, and thus in Eq. 4.30 one can replace 
$(t - s)$ with $s$.
\\
b) For arbitrary charging times and short discharging times, the solutions differ in that instead of $(t - s)^p$ the expression $s^p$ is used in the approximate solution.  However, as long as $p$ is very small and $(t - s)^p$ is slowly changing, no significant error will be caused by this.\\
If $(t - s)$ is replaced by $s$, the sum can be easily evaluated and Eq. 4.29 is obtained.
With regard to the approximation, what has been said under b) applies.
\\
Castro thus concludes that the Gross solution is generally usable and then has the advantage of great simplicity. The conclusions drawn from it, especially about the behavior with very short discharge times, remain strictly.
Indeed, later in 1940,  Gross \cite{Gross 1940} outlined the fact that De Oliveira 
Castro  has provided a more rigorous solution of the problem dealt by him
 in an approximate way in 1937 \cite{Gross 1937}.   

\subsection{The advent of  the Mittag-Leffler \\ functions
in linear viscoelasticity}
As we told at the beginning of this section,
Gross devoted  his attention  to the linear theory of viscoelastity based on the electro-mechanical analogy.
This was mostly since mid 1940's to mid 1950's, see  
\cite{Gross 1947,Gross 1948,Gross BOOK53}.
Being aware of the Mittag-Leffler function from the paper
\cite{Castro 1939},
 Gross  noted  in his 1947 paper
\cite{Gross 1947}  that this function is present in a viscoelastic model both in creep and relaxation (with different characteristic times) and provided  its spectral density  with the corresponding plots, as it was described in Section 2 and Fig 1.
\\
This fact has inspired Mainardi  in his PhD thesis carried out in the late 1960
at the University of Bologna under the supervision of Prof. Caputo.
As a matter of fact Mainardi was able to provide a plot of the function
$E_{\alpha}(-t^\alpha)$ for the first time in the literature, as one can see in the
1971  papers by Caputo and Mainardi 
\cite{Caputo-Mainardi 1971a,Caputo-Mainardi 1971b}
by introducing the so-called fractional Zener model.
More precisely in the survey \cite{Caputo-Mainardi 1971b}   families of 
 models of viscoelastic bodies were introduced generalizing the classical stress-strain relationship  of the mechanical models  by replacing ordinary derivatives with derivatives of non integer order.
For more details the reader is referred i.e. to Mainardi's book 
\cite{Mainardi BOOK10} published in 2010.
\\
We finally note that on the late 1960's  the  only  handbook 
(in English)  
dealing  with Mittag Leffler functions was that of  the  BATEMAN project
\cite{BATEMAN 55}, and  moreover  marginally in the chapter devoted  to "miscellaneous functions".
\section{Conclusions}
\label{S3} \vspace{-4pt}
\noindent
After a short survey of the basic properties of the Mittag-Leffler  functions 
 we have  shown the key role  of these  functions
 in dielectrical and mechanical processes, 
as outlined formerly by Harold Davis and Bernhard Gross.
We have discussed how  these two researchers have promoted
the Mittag-Leffler  functions  at their times when these functions were
practically unknown in applied sciences.
 In this survey we have illustrated  their contributions so
 it is worth to recall them  as pioneers of the applications of the Mittag-Leffler functions outside mathematics.\\
 {Moreover, we have outlined how the same results can be obtained making use of the Laplace Transform, instead of the method of successive approximations. This strategy has immediately led  to the solution and to the appearance of the Mittag-Leffler function.} 
 \\ \\
{\bf Acknowledgments:} 
The research activity of both authors
has been carried out in the framework of the activities of the National Group of Mathematical Physics (GNFM, INdAM), Italy.
{The authors would like to thank the 3  anonymous reviewers for their helpful and constructive comments.}
\\
We outline that  Section 3 is partly taken from the 1936  book by Davis 
including fig 8, whereas the   solution  of  Cole's circuit has been 
obtained by us using the  Laplace transform in terms of the Mittag-Leffler function. Furthermore, most of Section 4 has been adapted  from our   translation from German to English of the 1939 paper by Castro with the solutions obtained by us using Laplace transforms.  
\\

\end{document}